# On Pseudocontractions in Cyclic Maps

M. De la Sen

*Abstract-* **This paper discusses a more general contractive condition for a class of extended 2-cyclic self-mappings on the union of a finite number of subsets of a metric space which are allowed to have a finite number of successive images in the same subsets of its domain. If the space is uniformly convex and the subsets are non-empty, closed and convex then all the iterates converge to a unique closed limiting finite sequence which contains the best proximity points of adjacent subsets and reduces to a unique fixed point if all such subsets intersect.**

## I. INTRODUCTION

Strict pseudocontractive mappings and pseudocontractive mappings in the intermediate sense formulated in the framework of Hilbert spaces have received a certain attention in the last years concerning their convergence properties and the existence of fixed points (see, for instance, [1-4] and references therein). Results about the existence of a fixed point are discussed in those papers. On the other hand, important attention has been paid during the last decades to the study of the convergence properties of distances in cyclic contractive self-mappings on $p$ subsets $A_i \subset X$ of a metric space $(X,d)$, or a Banach space $(X,\|\ \|)$. The cyclic self-mappings under study have been of standard contractive or weakly contractive types and of Meir-Keeler type. The convergence of sequences to fixed points and best proximity points of the involved sets (see, for instance, [6-21] and references therein). It has to be noticed that every nonexpansive mapping is a 0-strict pseudocontraction and also that strict pseudontractions in the intermediate sense are asymptotically nonexpansive, [2]. The uniqueness of the best proximity points to which all the sequences of iterates converge is proven in [7] for the extension of the contractive principle for cyclic self- mappings in either uniformly convex Banach spaces (then being strictly convex and reflexive, [5]) or in reflexive Banach spaces, [14]. The $p$ subsets $A_i \subset X$ of the metric space $(X,d)$, or the Banach space $(X,\|\ \|)$, where the cyclic self-mappings are defined are supposed to be non-empty, convex and closed. If the involved subsets have nonempty intersections then all best proximity points are just a unique fixed point allocated in the intersection of all the subsets and framework can be simply given on complete metric spaces. The research in [7] is centred on the case of the 2-cyclic self-mapping being defined on the union of two subsets of the metric space. Those results are extended in [8] for Meir- Keeler cyclic contraction maps and, in general, for the self-mapping $T: \bigcup_{i \in \bar{p}} A_i \to \bigcup_{i \in \bar{p}} A_i$ be a $p(\geq 2)$-cyclic self-mapping being defined on any number of subsets of the metric space with $\bar{p} := \{1,2,...,p\}$.

Other recent researches which have been performed in the field of cyclic maps are related to the introduction and discussion of the so-called cyclic representation of a set $M$, decomposed as the union of a set of nonempty sets as $M = \bigcup_{i=1}^{m} M_i$, with respect to an operator $f : M \to M$, [15]. Subsequently, cyclic representations have been used in [16] to investigate operators from $M$ to $M$ which are cyclic $\phi\varphi$-contractions, where $\varphi : R_{0+} \to R_{0+}$ is a given comparison function, $M \subset X$ and $(X,d)$ is a metric space. The above cyclic representation has also been used in [17] to prove the existence of a fixed point for a self-mapping defined on a complete metric space which satisfies a cyclic weak $\varphi$-contraction. In [19], a characterization of best proximity points is studied for individual and pairs of non self-mappings $S,T : A \to B$, where $A$ and $B$ are nonempty subsets of a metric space. In general, best proximity points do not fulfil the usual "best proximity" condition $x = Sx = Tx$ under this framework. However, best proximity points are proven to jointly globally optimize the mappings from $x$ to the distances $d(x,Tx)$ and $d(x,Sx)$. Furthermore, a class of cyclic $\varphi$-contractions, which contain the cyclic contraction maps as a subclass, has been proposed in [19] in order to investigate the convergence and existence results of best proximity points in reflexive Banach spaces completing previous related results in [7]. Also, the existence and uniqueness of best proximity points of cyclic $\varphi$-contractive self-mappings in reflexive Banach spaces has been investigated in [20]. This paper is devoted to investigate the convergence properties and the existence of fixed points of a generalized version of pseudocontractve, strict pseudocontractive and asymptotically pseudocontractive in the intermediate sense in the more general framework of metric spaces. The case of 2-cyclic pseudocontractive self-mappings is also considered. The combination of constants defined the contraction may be different on each of the subsets and only the product of all the constants is requested to be less than unity. On the other hand, the self-mapping can perform a number of iterations on each of the subsets before transferring its image to the next adjacent subset of the 2-cyclic self-mapping. The existence of a unique closed finite limiting sequence on any sequence of iterates from any initial point in the union of the subsets is proven if $X$ is a uniformly convex Banach space and all the subsets of $X$ are nonempty, convex and closed. Such a limiting sequence is of size $q \geq p$ (with the inequality being strict if there is at least one iteration with image in the same subset as its domain) where $p$ of its elements (all of them if $q = p$) are best proximity points between adjacent subsets. In the case that all the subsets $A_i \subset X$ intersect, the above limit sequence reduces to a unique fixed point allocated within the intersection of all such subsets.

## II. ASYMPTOTIC CONTRACTIONS AND SEMICONTRACTIONS IN THE INTERMEDIATE SENSE

If $H$ is a real Hilbert space with an inner product $\langle .,. \rangle$ and a norm $\|.\|$ and $A$ is a nonempty closed convex subset of $H$ then the self-mapping $T : A \to A$ is said to be an asymptotically $\beta$-strictly pseudocontractive self-mapping in the intermediate sense for some $\beta \in [0,1)$ if

M. De la Sen is with the Institute of Research and Development of Processes, UPV/EHU, Bilbao, SPAIN (e-mail: manuel.delasen@ ehu.es).

$$\limsup_{n\to\infty} \sup_{x,y\in C} \left( \|T^n x - T^n y\|^2 - \alpha_n \|x-y\|^2 - \beta \|(I-T^n)x - (I-T^n)y\|^2 \right) \le 0$$

, $\forall x, y \in A$ (2.1)

for some sequence $\{\alpha_n\} \subset [1,\infty)$, $\alpha_n \to 1$ as $n \to \infty$, [1-5]. Such a concept was firstly introduced in [1]. If (2.1) holds for $\beta = 1$ then $T: A \to A$ is said to be an asymptotically pseudocontractive self-mapping in the intermediate sense. Finally, if $\alpha_n \to \alpha \in [0,1)$ as $n \to \infty$ then $T: A \to A$ is asymptotically $\beta$-strictly contractive in the intermediate sense, respectively, asymptotically contractive in the intermediate sense if $\beta = 1$. If (2.1) is changed to the stronger condition:

$$\left( \|T^n x - T^n y\|^2 - \alpha_n \|x-y\|^2 - \beta \|(I-T^n)x - (I-T^n)y\|^2 \right) \le 0;$$

$\forall x, y \in A$, $n \in N$ (2.2)

then the above concepts translate into $T: A \to A$ being an asymptotically $\beta$-strictly pseudocontractive self-mapping, an asymptotically pseudocontractive self-mapping and asymptotically contractive, respectively. Note that (2.1) is equivalent to:

$$\|T^n x - T^n y\|^2 \le \alpha_n \|x-y\|^2 + \beta \|(I-T^n)x - (I-T^n)y\|^2 + \xi_n;$$

$\forall x, y \in A$, $\forall n \in N$ (2.3)

or, equivalently,

$$\langle T^n x - T^n y, x-y \rangle \le \frac{1}{2}\left[ (k_n+1)\|x-y\|^2 + \xi_n \right]$$

; $\forall x, y \in A$, $n \in N$ (2.4)

where

$$\xi_n := \max\left\{ 0, \sup_{x,y\in C}\left( \|T^n x - T^n y\|^2 - \alpha_n \|x-y\|^2 - \beta \|(I-T^n)x - (I-T^n)y\|^2 \right) \right\}$$

; $\forall n \in N$ (2.5)

Note that the high-right-hand-side term $\|(I-T^n)x - (I-T^n)y\|^2$ of (2.1) is expanded as follows for any $x, y \in C$:

$$\|x-y\|^2 + \|T^n x - T^n y\|^2 - 2\|x-y\|\|T^n x - T^n y\|$$
$$\le \|(I-T^n)x - (I-T^n)y\|^2 = \langle x-T^n x, y-T^n y \rangle^2 = \langle x-y, T^n x - T^n y \rangle^2$$
$$= \|x-y\|^2 + \|T^n x - T^n y\|^2 + 2\langle T^n x - T^n y, x-y \rangle$$
$$= \langle x-y, T^n x - T^n y \rangle\langle x-y, T^n x - T^n y \rangle$$
$$\le \|x-y\|^2 + \|T^n x - T^n y\|^2 + 2\left|\langle T^n x - T^n y, x-y \rangle\right|$$
$$\le \|x-y\|^2 + \|T^n x - T^n y\|^2 + 2\|x-y\|\|T^n x - T^n y\| \quad (2.6)$$

The objective of this paper is to discuss the various pseudocontractive in the intermediate sense concepts in the framework of metric spaces and also to generalize them to the $\beta$- parameter to eventually be generalized to a sequence $\{\beta_n\}$ in $[0,1)$. Now, if instead of areal Hilbert space $H$ endowed with with an inner product $\langle .,.\rangle$ and a norm $\|.\|$, we deal with any generic Banach space $(X, \|.\|)$ then its norm induces an homogeneous and translation invariant metric $d: X \times X \to R_{0+}$ defined by $d(x,y) = d(x-y, 0) = \|x-y\|^{1/2}$; $\forall x, y \in A$ so that (2.6) takes the form:

$$d^2(x,y) + d^2(T^n x, T^n y) - 2d(x,y)d(T^n x, T^n y)$$
$$\le \|(I-T^n)x - (I-T^n)y\|^2 = d^2(x-y-(T^n x - T^n y), 0) = d^2(x-y, T^n x - T^n y)$$
$$\le (d(x-y,0) + d(T^n x - T^n y, 0))^2 = (d(x,y) + d(T^n x, T^n y))^2$$
$$= d^2(x,y) + d^2(T^n x, T^n y) + 2d(x,y)d(T^n x, T^n y); \forall x, y \in A$$
(2.7)

Define
$$\mu_n(x,y) := \min\left(\rho \in [-1,1]: d^2(x-y, T^n x - T^n y) \le d^2(x,y) + d^2(T^n x, T^n y) + 2\rho d(x,y)d(T^n x, T^n y)\right)$$
(2.8)

; $\forall x, y \in A$, $\forall n \in N$ which exists since it follows from (2.7)

$$\{1\} \subset \left\{\rho \in R: \|(I-T^n)x - (I-T^n)y\|^2 \le d^2(x,y) + d^2(T^n x, T^n y) + 2\rho d(x,y)d(T^n x, T^n y)\right\} (\ne \emptyset)$$
(2.9)

The following result holds related to the discussion (2.7)-(2.9) in metric spaces:

**Theorem 2.1.** Let $(X,d)$ be a metric space endowed with a homogeneous translation- invariant norm $d: X \times X \to R_{0+}$ and let $T: X \to X$ be a self-mapping. Assume that the constraint below holds:

$$d^2(T^n x, T^n y) \le \alpha_n(x,y)d^2(x,y) + \beta_n(x,y)\left(d^2(x,y) + d^2(T^n x, T^n y)\right)$$
$$+ 2\mu_n(x,y)\beta_n(x,y)d(x,y)d^2(T^n x, T^n y) + \xi_n(x,y)$$
; $\forall x, y \in X$, $\forall n \in N$ (2.10)

with
$$\xi_n = \xi_n(x,y)$$
$$:= \max\left(0, (1-\beta_n(x,y))d^2(T^n x, T^n y) - (\alpha_n(x,y)+\beta_n(x,y))d^2(x,y) - 2\mu_n(x,y)\beta_n(x,y)d(x,y)d^2(T^n x, T^n y)\right)$$
$$\to 0 \quad ; \forall x, y \in X \text{ as } n \to \infty \quad (2.11)$$

for some parameterizing real sequences $\alpha_n = \alpha_n(x,y)$, $\beta_n = \beta_n(x,y)$ and $\mu_n = \mu_n(x,y)$ satisfying for any $n \in N$:

$$\{\alpha_n(x,y)\} \subset [0,\infty), \quad \{\mu_n(x,y): x \ne y\} \subset \left(-\infty, \frac{1-\beta_n(x,y)}{2\beta_n(x,y)d(x,y)}\right],$$

$\mu_n(0,0) = 0$, $\{\beta_n(x,y)\} \subset [0,1)$  ; $\forall x, y \in X$, $\forall n \in N$
(2.12)

Then, the following properties hold:

**(i)** $\exists \lim_{n\to\infty} d(T^n x, T^n y) \le d(x,y)$ for any $x, y \in X$ satisfying the conditions:

either $\dfrac{\alpha_n(x,y) + \beta_n(x,y)}{1-\beta_n(x,y)(1+2\mu_n(x,y)d(x,y))} \le 1$ ; $\forall n \in N$, or

$\dfrac{\alpha_n(x,y) + \beta_n(x,y)}{1-\beta_n(x,y)(1+2\mu_n(x,y)d(x,y))} \ge 1$; $\forall x, y \in X$, $\forall n \in N$ (2.13)

$\alpha_n(x,y) + 2\beta_n(x,y)(1+\mu_n(x,y)d(x,y)) \to 1$ as $n \to \infty$ ; $\forall x, y \in X$ (2.14)

and, if (2.13) holds for $\forall x(=y) \in X$, and, furthermore:

$\alpha_n(x,Tx) + 2\beta_n(x,Tx)(1+\mu_n(x,Tx)d(x,Tx)) \to 1$ as $n \to \infty$; $\forall x \in X$ (2.15)

then $\exists \lim_{n\to\infty} d(T^n x, T^{n+1} x) \le d(x,Tx)$. If (2.13)-(2.14) hold for any $x, y \in X$ then $T: X \to X$ is asymptotically nonexpansive.

**(ii)** If, in addition, $(X,d)$ is a convex metric space then $\lim_{n\to\infty} d(T^n x, T^{n+m+1} x) = 0$; $\forall x \in X$, $\forall m \in \mathbf{N}$ and $\exists \lim_{n\to\infty} d(T^{n+m} x, z)$; $\forall x, y \in X$, $\forall m \in \mathbf{N}$. Furthermore, if $(X,d)$ is complete then $\{T^n x\}_{n \in \mathbf{N}}$ converges in $X$. If, furthermore, $T: X \to X$ is continuous then each sequence $\{T^n x\}_{n \in \mathbf{N}}$ converges to a fixed point in $X$.

*Proof*: Note that (2.9) is well-posed since the minimum always exists since the set over which such a minimum is computed contains $\{1\}$. The constraint (2.12) leads to:

$$d^2(T^n x, T^n y) \leq \frac{\alpha_n(x,y) + \beta_n(x,y)}{1 - \beta_n(x,y)(1 + 2\mu_n(x,y)d(x,y))} d^2(x,y) + \frac{\xi_n(x,y)}{1 - \beta_n(x,y)(1 + 2\mu_n(x,y)d(x,y))}$$

(2.16)

If (2.11)-(2.14) hold for some $x, y \in X$ then $\exists \lim_{n\to\infty} d(T^n x, T^n y) \leq d(x,y)$ from (2.16).

If the constraints (2.11)-(2.14) hold for all $x, y \in X$ then $\exists \lim_{n\to\infty} d(T^n x, T^n y) = L = L(x,y) \leq d(x,y)$. Assume that $(X,d)$ is a convex metric space. Then, the triangle inequality of distances reduces to an identity for any pairs of elements in $X$. Then, from the above limit property of distances and proceeding by complete induction by assuming that $\exists \lim_{n\to\infty} d(T^n x, T^{n+j} x) = 0$; $\forall x \in X$, $\forall j (\leq m) \in \mathbf{N}$ for any given $m \in \mathbf{N}$, then one gets by the triangle equality:

$0 \leftarrow d(T^n x, T^{n+m} x) = d(T^n x, T^{n+m+1} x) + d(T^{n+m} x, T^{n+m+1} x) \to d(T^n x, T^{n+m+1} x)$ as $n \to \infty$; $\forall x \in X \Rightarrow \exists \lim_{n\to\infty} d(T^n x, T^{n+m+1} x) = 0$; $\forall x \in X$ (2.17)

and

$d(T^n x, y) = d(T^n x, T^{n+m} x) + d(T^{n+m} x, y) \to d(T^{n+m} x, y) = L_1(x,y)$
; $\forall x, y \in X$, $\forall m \in \mathbf{N}$ (2.18)

so that $\exists \lim_{n\to\infty} d(T^{n+m} x, y)$; $\forall x, y \in X$, $\forall m \in \mathbf{N}$ and $T: X \to X$ is a asymptotically nonexpansive. Property (i) has been proven. Property (ii) is proven as follows. Assume that $\{T^n x\}$ does not converge for some $x \in X$. Thus, $\neg \exists x, z \in X$ such that $T^n x \to z$ as $n \to \infty$ (this includes that there is no $x$ in $X$ such that $x = T^n x$ for some $n \in \mathbf{N}$). Since the metric space is convex then there is $z_1 \in X$ such that $d(T^n x, T^{n+1} x) = d(T^n x, z_1) + d(T^{n+1} x, z_1) \to 0$ as $n \to \infty$ and the contradiction $T^n x \to z$ as $n \to \infty$ follows. Then, $\{T^n x\}$ has a limit in $X$ for any $x \in X$ since $(X,d)$ is complete. If, in addition, $T: X \to X$ is continuous then each sequence $\{T^n x\}$ converges to a fixed point in $X$ since $z \leftarrow T^n x \leftarrow T^{n+1} x = T(T^n x) \to Tz$, hence $z = Tz$ is a fixed point of $T: X \to X$. □

The following result extends Theorem 2.1 for a modification of the contractive condition (2.10):

**Theorem 2.2**. Let $(X,d)$ be a metric space endowed with a homogeneous translation-invariant norm $d: X \times X \to \mathbf{R}_{0+}$ and let $T: X \to X$ be a self-mapping. Assume that the constraint below holds:

$d^2(T^n x, T^n y) \leq \alpha_n(x,y) d^2(x,y) + \beta_n(x,y)(d^2(x,y) + d^2(T^n x, T^n y))$
$+ 2\mu_n(x,y)\beta_n(x,y)d(x,y)d(T^n x, T^n y) + \xi_n(x,y)$; $\forall x, y \in X$, $\forall n \in \mathbf{N}$ (2.19)

with

$\xi_n = \xi_n(x,y)$
$:= \max(0, (1 - \beta_n(x,y))d^2(T^n x, T^n y) - (\alpha_n(x,y) + \beta_n(x,y))d^2(x,y) - 2\mu_n(x,y)\beta_n(x,y)d(x,y)d^2(T^n x, T^n y))$
$\to 0$; $\forall x, y \in X$ as $n \to \infty$ (2.20)

for some parameterizing real sequences $\alpha_n = \alpha_n(x,y)$, $\beta_n = \beta_n(x,y)$ and $\mu_n = \mu_n(x,y)$ satisfying for any $n \in \mathbf{N}$:

$\{\alpha_n(x,y)\} \subset [0, \infty)$,

$\{\mu_n(x,y): x \neq y\} \subset \left(-\infty, \frac{1 - \beta_n(x,y)}{2\beta_n(x,y)}\right]$, $\{\beta_n(x,y)\} \subset [0,1)$ (2.21)

; $\forall x, y \in X$, $\forall n \in \mathbf{N}$. Then the following properties hold:

**(i)** $\exists \lim_{n\to\infty} d(T^n x, T^n y) \leq d(x,y)$ for any $x, y \in X$ satisfying the conditions:

Either $\frac{\alpha_n(x,y) + \beta_n(x,y)}{1 - \beta_n(x,y)(1 + 2\mu_n(x,y)d(x,y))} \leq 1$; $\forall n \in \mathbf{N}$, or

$\frac{\alpha_n(x,y) + \beta_n(x,y)}{1 - \beta_n(x,y)(1 + 2\mu_n(x,y)d(x,y))} \geq 1$ (2.22)

$\alpha_n(x,y) + 2\beta_n(x,y)(1 + \mu_n(x,y)) \to 1$; $\forall x, y \in X$ as $n \to \infty$ (2.23)

and, if (2.22) holds for all $x, y (= x) \in X$, and, furthermore:

$\alpha_n(x, Tx) + 2\beta_n(x, Tx)(1 + \mu_n(x, Tx)) \to 1$; $\forall x \in X$ as $n \to \infty$ (2.24)

then $\exists \lim_{n\to\infty} d(T^n x, T^{n+1} x) \leq d(x, Tx)$.

If (2.23) holds for any $x, y \in X$ then $T: X \to X$ is asymptotically nonexpansive.

**(ii)** If, in addition $(X,d)$ is a convex complete metric space and $T: X \to X$ is continuous then $T: X \to X$ has a fixed point in $X$.

*Proof*: Four cases can occur implied by (2.19) for each $x, y \in X$ and $n \in \mathbf{N}$, namely:

a) $\left(d(T^n x, T^n y) \geq d(x,y)\right) \wedge \left(\frac{1 - \beta_n(x,y)}{2\beta_n(x,y)} \geq \mu_n(x,y) \geq 0\right)$. Then, (2.19) implies that:

$d^2(T^n x, T^n y) \leq \frac{\alpha_n(x,y) + \beta_n(x,y)}{1 - \beta_n(x,y)(1 + 2\mu_n(x,y))} d^2(x,y) + \frac{\xi_n(x,y)}{1 - \beta_n(x,y)(1 + 2\mu_n(x,y))}$

(2.25)

b) $\left(d(T^n x, T^n y) < d(x,y)\right) \wedge \left(\frac{1 - \beta_n(x,y)}{2\beta_n(x,y)} \geq \mu_n(x,y) \geq 0\right)$. Then, (2.19) implies that:

$d^2(T^n x, T^n y) \leq \frac{\alpha_n(x,y) + \beta_n(x,y)(1 + 2\mu_n(x,y))}{1 - \beta_n(x,y)} d^2(x,y) + \frac{\xi_n(x,y)}{1 - \beta_n(x,y)}$ (2.26)

c) $\left(d(T^n x, T^n y) \geq d(x,y)\right) \wedge \left((\beta_n(x,y) - 1)/2\beta_n(x,y) < \mu_n(x,y) < 0\right)$. Then, (2.19) implies that:

$$d^2(T^n x, T^n y) \le \frac{\alpha_n(x,y) + \beta_n(x,y)(1 - 2|\mu_n(x,y)|)}{1 - \beta_n(x,y)} d^2(x,y)$$
$$+ \frac{\xi_n(x,y)}{1 - \beta_n(x,y)} \quad (2.27)$$

d)
$$(d(T^n x, T^n y) < d(x,y)) \wedge ((\beta_n(x,y)-1)/2\beta_n(x,y) < \mu_n(x,y) < 0).$$

Then, (2.19) implies that:

$$d^2(T^n x, T^n y) \le \frac{\alpha_n(x,y) + \beta_n(x,y)}{1 - \beta_n(x,y)(1 - 2|\mu_n(x,y)|)} d^2(x,y)$$
$$+ \frac{\xi_n(x,y)}{1 - \beta_n(x,y)(1 - 2|\mu_n(x,y)|)} \quad (2.28)$$

If the limiting condition (2.23) holds for any $x, y \in X$, subject to (2.23)-(2.24), leads from any constraints (2.25)-(2.28) to:

$$\exists \lim_{n \to \infty} d(T^n x, T^n y) \le d(x,y)$$
$$\Rightarrow \lim_{n \to \infty} (d^2(T^{n+1} x, T^{n+1} y) - d^2(T^n x, T^n y)) = 0 \,; \forall x, y \in X \quad (2.29)$$

so that $T: X \to X$ is asymptotically nonexpansive. Using the same proving arguments as in Theorem 2.1, one proves that (2.29) implies that $d^2(T^n x, T^n y) \to 0$ as $n \to \infty$; $\forall x, y \in X$ and if the metric space $(X,d)$ is convex and complete then $T: X \to X$ has a fixed point in $X$ to which all sequences $\{T^n x\}$ converge; $\forall x \in X$. □

The concepts of pseudocontraction and asymptotic pseudocontraction in the intermediate sense motivated by (2.7)-(2.9) are revisited as follows in the context of metric spaces and then linked to Theorem 2.2:

**Definition 2.3**. Assume that $(X,d)$ is a complete metric space with $d: X \times X \to \boldsymbol{R}_{0+}$ being a homogeneous translation-invariant metric. Thus, $T: A \to A$ is asymptotically $\beta$ - strictly pseudocontractive in the intermediate sense if:

$$\limsup_{n \to \infty} ((1 + 2\mu_n - \beta_n)d^2(T^n x, T^n y) - (\alpha_n + \beta_n)d^2(x,y)) \le 0$$
$$; \forall x, y \in A \quad (2.30)$$

for some real sequences $\{\mu_n\}$, $\{\alpha_n\}$ and $\{\beta_n\}$ satisfying:
$\{\mu_n\} \subset [\mu, \infty), \{\alpha_n\} \subset [1, \infty), \{\beta_n\} \subset [0, 1 + 2\mu_n], \mu \in \boldsymbol{R}_{0+}$,
$\alpha_n \to \alpha = 1, \beta_n \to \beta \in [0,1)$ as $n \to \infty$; $\forall x, y \in A$, $\forall n \in \boldsymbol{N}$
(2.31)
□

**Definition 2.4**. $T: A \to A$ is asymptotically pseudocontractive in the intermediate sense if $\beta = 1$ in Definition 2.3. □

**Definition 2.5**. $T: A \to A$ is asymptotically $\beta$ - strictly contractive in the intermediate sense if $\alpha_n \in [0,1)$; $\ni n \in \boldsymbol{N}$, $\alpha_n \to \alpha \in [0,1)$ as $n \to \infty$ and $\beta \in [0,1)$ in Definition 2.3. □

**Definition 2.5**. $T: A \to A$ is asymptotically contractive in the intermediate sense if $\alpha \in [0,1)$ and $\beta = 1$ in Definition 2.4. □

The following result being supported by Theorem 2.2 relies on the concepts of asymptotically contractive and pseudocontractive self-mappings in the intermediate sense. Therefore, it is assumed that $\{\alpha_n(x,y)\} \subset [1, \infty)$.

**Theorem 2.6**. Let $(X,d)$ be a convex complete metric space endowed with a homogeneous translation-invariant norm $d: X \times X \to \boldsymbol{R}_{0+}$ and let and let $T: A \to A$ be a self-mapping where $A$ is a nonempty, closed and convex subset of $X$. Assume that the constraint below holds:

$$d^2(T^n x, T^n y) \le \alpha_n(x,y)d^2(x,y) + \beta_n(x,y)(d^2(x,y) + d^2(T^n x, T^n y))$$
$$+ 2\mu_n(x,y)\beta_n(x,y)d(x,y)d(T^n x, T^n y) + \xi_n(x,y); \forall x, y \in A, \forall n \in \boldsymbol{N} \quad (2.32)$$

with
$\xi_n = \xi_n(x,y)$
$:= \max(0, (1 - \beta_n(x,y))d^2(T^n x, T^n y) - (\alpha_n(x,y) + \beta_n(x,y))d^2(x,y) - 2\mu_n(x,y)\beta_n(x,y)d(x,y)d^2(T^n x, T^n y))$
$\to 0$ ; $\forall x, y \in A$ as $n \to \infty$ (2.33)

for some parameterizing real sequences $\alpha_n = \alpha_n(x,y)$, $\beta_n = \beta_n(x,y)$ and $\mu_n = \mu_n(x,y)$ satisfying for any $n \in \boldsymbol{N}$:
$\{\alpha_n(x,y)\} \subset [1, \infty)$,
$\{\mu_n(x,y): x \ne y\} \subset \left(-\infty, \frac{1 - \beta_n(x,y)}{2\beta_n(x,y)}\right], \{\beta_n(x,y)\} \subset [0, \beta) \subset [0,1]$ (2.34)
; $\forall x, y \in X$, $\forall n \in \boldsymbol{N}$. Then, $\exists \lim_{n \to \infty} d(T^n x, T^n y) \le d(x,y)$ for any $x, y \in X$ satisfying the conditions:

$$\frac{\alpha_n(x,y) + \beta_n(x,y)}{1 - \beta_n(x,y)(1 + 2\mu_n(x,y)d(x,y))} \ge 1\,;$$
$\alpha_n(x,y) + 2\beta_n(x,y)(1 + \mu_n(x,y)) \to 1$; $\forall x, y \in X$ as $n \to \infty$ (2.35)

Then, the following properties hold:

(i) $T: A \to A$ is asymptotically $\beta$-strictly pseudocontractive in the intermediate sense if Eqs. (2.34) hold with $0 < \beta < 1$, $1 \le \alpha_n \to 1$ as $n \to \infty$ and (2.32) holds for $\mu_n \to -1$ as $n \to \infty$. Also, $T: A \to A$ has a fixed point in $A$ if $T: X \to X$ is continuous.

(ii) If $\beta = 1$, $1 \le \alpha_n \to 1$ as $n \to \infty$ and (2.32) holds for $\mu_n \to -1$ as $n \to \infty$ then $T: A \to A$ is asymptotically pseudocontractive in the intermediate sense. Also, $T: A \to A$ has a fixed point in $A$ if $T: X \to X$ is continuous.

(iii) If (2.34) is modified as $0 \le \alpha_n(n \in \boldsymbol{N}) \to \alpha \in [0,1)$ and $\mu_n \to -1$ as $n \to \infty$, and $\beta \in [0,1)$ then $T: A \to A$ is asymptotically $\beta$ - strictly contractive in the intermediate. Also, $T: A \to A$ has a fixed point in $A$ if $T: X \to X$ is continuous.

(iv) If (2.34) is modified as $0 \le \alpha_n(n \in \boldsymbol{N}) \to \alpha \in [0,1)$ and $\mu_n \to -1$ as $n \to \infty$, and $\beta = 1$ then $T: A \to A$ is asymptotically contractive. Also, $T: A \to A$ has a fixed point in $A$ if $T: X \to X$ is continuous.

*Proof*: It follows from Definitions 2.2-2.5 and the fact that Theorem 2.2 holds with the condition (2.32) on $A$ with $\alpha_n \to 1$ and $\mu_n \to -1$ as $n \to \infty$ subject to (2.27) or (2.28) for each $x, y \in A$, $n \in \boldsymbol{N}$. □

For uniqueness of the fixed point in the various results of this section, we state the subsequent result being a Corollary to Theorems 2.1, 2.2 and 2.6:

**Corollary 2.7**. If $A$ is a nonempty closed convex subset of $X$, $T : A \to A$ containing $\{0\}$ and $(X,d)$ is a metric space then there is a unique fixed point in $A$ under the conditions of existence of at least a fixed point given in Theorems 2.1. 2.2 and 2.6.

*Proof*: Assume that $x = Tx, y = Ty(\neq x)$ are two fixed points in $A$. Since $A$ is nonempty, closed and convex, $\exists z \in A$ such that $z \in (x,y) \subset A$, and
$$d(Tx,Ty) = 2d(x,z) = 2d(y,z) = 2d(Tx,z) = 2d(Ty,z) \quad (2.36)$$

implying that $d(x,z) = d(y,z)$ and, since the metric is translation-invariant, then $d(x,y) = 0$ so that $x = y$ leading to a contradiction. Hence, the corollary. □

## II. ASYMPTOTIC PROPERTIES IN THE INTERMEDIATE SENSE OF CYCLIC SELF-MAPPINGS

Let $A, B \subset X$ be nonempty subsets of $X$, $T : A \cup B \to A \cup B$ is cyclic self-mapping if $T(A) \subseteq B$ and $T(B) \subseteq A$. Assume that contractive condition (2.19) is modified as follows:
$$d^2(T^n x, T^n y) \leq \alpha_n(x,y) d^2(x,y) + \beta_n(x,y)\left(d^2(x,y) + d^2(T^n x, T^n y)\right)$$
$$+ 2\mu_n(x,y)\beta_n(x,y) d(x,y) d(T^n x, T^n y) + \xi_n(x,y) + \gamma_n(x,y) D^2$$
$$; \forall x \in A, y \in B, \forall n \in N \quad (3.1)$$

where $\{\gamma_n(x,y)\} \in [0, \infty)$ and $D = dist(A, B) \geq 0$. If $A \cap B \neq \emptyset$ then $D = 0$ and Theorem 2.2 holds with the replacement $A \to A \cap B$. Then, if $A$ and $B$ are closed and convex then there is a unique fixed point of $T : A \cup B \to A \cup B$ in $A \cap B$. In the following, we consider the case that $A \cap B = \emptyset$ so that $D > 0$. The following result based on Theorem 2.6 holds:

**Theorem 3.1**. Let $(X, \|\ \|)$ be a reflexive Banach space with a norm-induced homogeneous translation-invariant norm $d : X \times X \to \mathbf{R}_{0+}$, where $A$ and $B \subset X$ are nonempty, closed and convex subsets of $X$ such, $T : A \cup B \to A \cup B$ is cyclic self-mapping if $T(A) \subseteq B$ and $T(B) \subseteq A$. Define the sequence $\{k_n\}_{n \in N} \subset [0, \infty)$ as follows:

$$k_n = k_n(x,y) := \begin{cases} \frac{\alpha_n(x,y) + \beta_n(x,y)}{1 - \beta_n(x,y)(1 + 2\mu_n(x,y))} & \text{if } \frac{1 - \beta_n(x,y)}{2\beta_n(x,y)} \geq \mu_n(x,y) \geq 0 \\ \frac{\alpha_n(x,y) + \beta_n(x,y)(1 + 2\mu_n(x,y))}{1 - \beta_n(x,y)} & \text{if } (\beta_n(x,y) - 1)/2\beta_n(x,y) < \mu_n(x,y) < 0 \end{cases}$$
if $d(T^n x, T^n y) \geq d(x,y)$ (3.2)
and
$$k_n = k_n(x,y) := \begin{cases} \frac{\alpha_n(x,y) + \beta_n(x,y)(1 + 2\mu_n(x,y))}{1 - \beta_n(x,y)} & \text{if } \frac{1 - \beta_n(x,y)}{2\beta_n(x,y)} \geq \mu_n(x,y) \geq 0 \\ \frac{\alpha_n(x,y) + \beta_n(x,y)}{1 - \beta_n(x,y)(1 + 2\mu_n(x,y))} & \text{if } (\beta_n(x,y) - 1)/2\beta_n(x,y) < \mu_n(x,y) < 0 \end{cases}$$
if $d(T^n x, T^n y) < d(x,y)$ (3.3)
; $\forall x \in A, y \in B$ for some parameterizing real sequences $\alpha_n = \alpha_n(x,y), \beta_n = \beta_n(x,y)$ and $\mu_n = \mu_n(x,y)$, $\xi_n = \xi_n(x,y)$ satisfying for any $n \in N$:
$\{\alpha_n(x,y) := 1 - k_n(x,y)\} \subset [0, \infty)$,
$$\{\mu_n(x,y): x \neq y\} \subset \left(-\infty, \frac{1 - \beta_n(x,y)}{2\beta_n(x,y)}\right] \quad (3.4)$$
$\{\beta_n(x,y)\} \subset [0, \beta) \subset [0,1], \{\gamma_n(x,y)\} \in [0, \infty) \quad (3.5)$

$\xi_n = \xi_n(x,y)$
$:= \max(0, (1 - \beta_n(x,y)) d^2(T^n x, T^n y) - (\alpha_n(x,y) + \beta_n(x,y)) d^2(x,y) - 2\mu_n(x,y)\beta_n(x,y) d(x,y) d^2(T^n x, T^n y))$
$\to 0 \quad ; \forall x \in A, y \in B \text{ as } n \to \infty \quad (3.6)$
Then, the following properties hold:

**(i)** $T : A \cup B \to A \cup B$ is asymptotically $\beta$-strictly pseudocontractive in the intermediate sense if (3.5) holds with $0 < \beta < 1$, $1 \leq \alpha_n \to 1$ as $n \to \infty$ and (3.6) holds for $\mu_n \to -1$ as $n \to \infty$. Also, $T : A \cup B \to A \cup B$ has a best proximity point in $A$ and a best proximity point in $B$ to which the sequences $\{T^{2n} x\}$ and $\{T^{2n+1} x\}$ converge if $T : A \cup B \to A \cup B$ is continuous.

**(ii)** If $\beta = 1$, $1 \leq \alpha_n \to 1$ as $n \to \infty$ and (3.6) holds for $\mu_n \to -1$ as $n \to \infty$ then $T : A \cup B \to A \cup B$ is asymptotically pseudocontractive in the intermediate sense. Also, $T : A \cup B \to A \cup B$ has a best proximity point in $A$ and a best proximity point in $B$ to which the sequences $\{T^{2n} x\}$ and $\{T^{2n+1} x\}$ converge if $T : A \cup B \to A \cup B$ is continuous.

**(iii)** If (3.4) is modified as $0 \leq \alpha_n (n \in N) \to \alpha \in [0,1)$ and $\mu_n \to -1$ as $n \to \infty$, and $\beta \in [0,1)$ in (3.5) then $T : A \cup B \to A \cup B$ is asymptotically $\beta$-strictly contractive in the intermediate. Also, $T : A \cup B \to A \cup B$ has a best proximity point in $A$ and a best proximity point in $B$ to which the sequences $\{T^{2n} x\}$ and $\{T^{2n+1} x\}$ converge if $T : A \cup B \to A \cup B$ is continuous.

**(iv)** If (3.4) is modified as $0 \leq \alpha_n (n \in N) \to \alpha \in [0,1)$ and $\mu_n \to -1$ as $n \to \infty$, and $\beta = 1$ in (3.5) then $T : A \cup B \to A \cup B$ is asymptotically contractive. Also, $T : A \cup B \to A \cup B$ has a best proximity point in $A$ and a best proximity point in $B$ to which the sequences $\{T^{2n} x\}$ and $\{T^{2n+1} x\}$ converge if $T : A \cup B \to A \cup B$ is continuous.

*Proof*: If $T : A \cup B \to A \cup B$ is asymptotically $\beta$-strictly pseudocontractive in the intermediate sense under (3.5)-(3.7) with $0 < \beta < 1$, $1 \leq \alpha_n \to 1$ as $n \to \infty$ and (3.6) holds for $\mu_n \to -1$ as then $D \leq d(T^n x, T^n y)(\leq d(x,y) = D) \to D$ as $n \to \infty$ for any $x \in A$ and $y \in B$. Thus, for any $\varepsilon \in \mathbf{R}_+$, $\exists n_0 = n_0(\varepsilon) \in \mathbf{Z}_+$ such that $D \leq d(T^{2n} x, T^{2n} y) \leq D + \varepsilon$; $\forall n \geq n_0$ with $\{T^{2j} x\} \subset A$ and $\{T^{2j} y\} \subset B$; $\forall j \in \mathbf{Z}_+$ if $x \in A$ and $y \in B$ since $T(A) \subseteq B$ and $T(B) \subseteq A$. Furthermore, $\{T^{2n} x\} \subset \hat{A}$ and $\{T^{2n} x\} \subset \hat{B}$; $\forall n \geq n_0$ with $\hat{A}$ and $\hat{B}$ being nonempty closed and convex subsets of $A$ and $B$, respectively, which exists and are bounded since $D \leq d(T^n x, T^n y) \to D$ as $n \to \infty$ for any $x \in A$ and $y \in B$. Then, $D \leq d(T^{2n} x, T^{2n} y)(\leq d(x,y)) \to D = d(z, \omega)$ for some $z \in \hat{A}$ and $\omega \in \hat{B}$ since $(X, \|\ \|)$ is a reflexive Banach space (i.e. $(X,d)$ is a complete metric space under the nor-induced distance) and since $\hat{A}$ is nonempty, bounded, closed and convex and $\hat{B}$ is nonempty, closed and convex, so that $D = d(z, \omega)$ for some $z \in \hat{A}$ and $\omega \in \hat{B}$, [7]. Hence:

$$d(T^{2n}x, T^{2n}y) - d(T^{2n+2}x, T^{2n+2}y) = d(T^{2n}x, T^{2n}y) - d(T^{2n+2}x + (T^{2n+2}x - T^{2n}x), T^{2n}y + (T^{2n+2}y - T^{2n}y))$$
$$\to D - D = 0 \text{ as } n \to \infty \quad (3.7)$$

and one gets by taking $y = Tx$

$$d(T^{2n}x, T^{2n+1}x) \to d(T^{2n}x + (T^{2n+2}x - T^{2n}x), T^{2n}y + (T^{2n+3}x - T^{2n+1}x))$$
$$\text{as } n \to \infty \quad (3.8)$$

Thus, $\{T^{2n}x\}$ and $\{T^{2n+1}x\}$ are Cauchy sequences which converge, respectively, to $z_1 \in A$ and $z_2 \in B$ since $A$ and $B$ are closed sets since $d(z_1, z_2) = D$. If $T: A \cup B \to A \cup B$ is continuous then $z_2 = Tz_1$ are best proximity points of $T: A \cup B \to A \cup B$ since $z_2 \leftarrow T^{2n+1}x = T(T^{2n}x) \to Tz_1$ as $n \to \infty$. Property (i) has been proven. The remaining properties follow directly under similar arguments. □

*Remark 3.2.* Note that the existence of Theorem 3.1 of $x \in A$ and $y \in B$ such that $D = d(x, y)$ guaranteed if $A$ is nonempty bounded, closed and convex and $B$ is nonempty closed and convex is also guaranteed if $A$ is compact and $B$ is approximatively compact with respect to $A$, i.e. if every sequence $\{x_n\} \subset B$ such that $d(y, x_n) \to d(y, B)$ for $y \in A$ has a convergent subsequence, [7]. □

## ACKNOWLEDGMENTS

The author is grateful to the Spanish Ministry of Education by its partial support of this work through Grant DPI2009-07197. He is also grateful to the Basque Government by its support through Grants IT378-10 and SAIOTEK S-PE09UN12.

## REFERENCES

[1] D.R. Sahu, H.K. Xu and J.C. Yao, "Asymptotically strict pseudocontractive mappings in the intermediate sense", *Nonlinear Analysis, Theory, Methods and Applications*, Vol. 70, pp. 3502-3511, 2009.

[2] X. Qin, J.K. Kim and T. Wang, "On the convergence of implicit iterative processes for asymptotically pseudocontractive mappings in the intermediate sense", *Abstract and Applied analysis*, Vol. 2011, Article Number 468716, 18 pages, 2011, doi:10.1155/2011/468716.

[3] L.C. Ceng, A. Petrusel and J.C. Yao, "Iterative approximation of fixed points for asymptotically strict pseudocontarctive type mappings in the intermediate sense", *Taiwanese Journal of Mathematics*, Vol. 15, No. 2, pp. 587-606, 2011.

[4] P. Duan and J. Zhao, "Strong convergence theorems for system of equilibrium problems and asymptotically strict pseudocontractions in the intermediate sense", *Fixed Point Theory and Applications*, Vol. 2011, doi:10.1186/1687-1812-2011-13.

[5] X. Qin, S.M. Kang and R. P. Agarwal, "On the convergence of an implicit iterative process for generalized asymptotically quasi-nonexpansive mappings", *Fixed Point Theory and Applications*, Vol. 2010, Article Number 714860, 19 pages, doi:10.1155/2010/714860.

[6] W.A. Kirk, P.S. Srinivasan and P. Veeramani, "Fixed points for mappings satisfying cyclical contractive conditions", *Fixed Point Theory*, Vol. 4, No. 1, pp. 79-89, 2003.

[7] A.A. Eldred and P. Veeramani, "Existence and convergence of best proximity points", *J. Math. Anal. Appl.*, Vol 323, pp. 1001-1006, 2006.

[8] S. Karpagam and S. Agrawal, "Best proximity point theorems for p-cyclic Meir-Keeler contractions", *Fixed Point Theory and Applications*, Vol. 2009, Article Number 197308, 9 pages, 2009, doi:10.1155/2009/197308.

[9] C. Di Bari, T. Suzuki and C. Vetro, "Best proximity points for cyclic Meir-Keeler contractions", Nonlinear Analysis: Theory, Methods & Applications, Vol. 69, No. 11, pp. 3790-3794, 2008.

[10] M. De la Sen, "Linking contractive self-mappings and cyclic Meir-Keeler contractions with Kannan self-mappings", *Fixed Point Theory and Applications*, Vol. 2010, Article Number 572057, 23 pages, 2010, doi:10.1155/2010/572057.

[11] M. De la Sen, "Some combined relations between contractive mappings, Kannan mappings reasonable expansive mappings and T-stability", *Fixed Point Theory and Applications*, Vol. 2009, Article Number 815637, 25 pages, 2009, doi:10.1155/2009/815637.

[12] T. Suzuki, "Some notes on Meir-Keeler contractions and L-functions", *Bulletin of the Kyushu Institute of Technology,* No.53, pp. 12-13, 2006.

[13] M. Derafshpour, S. Rezapour and N. Shahzad, "On the existence of best proximity points of cyclic contractions", *Advances in Dynamical Systems and Applications*, Vol. 6, No. 1, pp. 33-40, 2011.

[14] Sh. Rezapour, M. Derafshpour and N. Shahzad, "Best proximity points of cyclic $\varphi$-contractions on reflexive Banach spaces", *Fixed Point Theory and Applications*, Vol. 2010 (2010), Article ID 946178, 7 pages, doi:10.1155/2010/046178.

[15] M.A. Al-Thagafi and N. Shahzad, "Convergence and existence results for best proximity points", *Nonlinear Analysis, Theory, Methods & Applications*, Vol. 70, No. 10, pp. 3665-3671, 2009.

[16] I.A. Rus, "Cyclic representations and fixed points", *Ann. T. Popoviciu Seminar Funct. Eq. Approx. Convexity*, Vol. 3, pp. 171-178. 2005.

[17] M. Pacurar and I.A. Rus, "Fixed point theory for cyclic $\varphi$-contraction", *Nonlinear Analysis- Theory, Methods & Applications*, Vol. 72, No. 3-4, pp. 1181-1187, 2010.

[18] E. Karapinar, "Fixed point theory for cyclic weak $\phi$-contraction", *Applied Mathematics Letters*, Vol. 24, pp. 822-825, 2011.

[19] N. Shazhad, S. Sadiq Basha and R. Jeyaraj, "Common best proximity points: global optimal solutions", *J. Optim. Theory Appl*, Vol. 148, No.1, pp. 69-78, 2011.

[20] C. Vetro, "Best proximity points: convergence and existence theorems for p-cyclic mappings", *Nonlinear Analysis- Theory, Methods & Applications*, Vol. 73, No. 7, pp. 2283-2291, 2010.